\documentclass{article}%
\usepackage{makeidx}
\usepackage{amsfonts}
\usepackage{amsmath}%
\usepackage{amssymb}%
\usepackage{graphicx}
\newtheorem{theorem}{Theorem}

\newenvironment{proof}[1][Proof]{\noindent\textbf{#1.} }{\ \rule{0.5em}{0.5em}}
\begin{document}

\title{Affine and Projective Universal Geometry}
\author{N J Wildberger\\School of Mathematics and Statistics\\UNSW Sydney 2052 Australia}
\maketitle
\date{}

\begin{abstract}
By recasting metrical geometry in a purely algebraic setting, both Euclidean
and non-Euclidean geometries can be studied over a general field with an
arbitrary quadratic form. Both an affine and a projective version of this new
theory are introduced here, and the main formulas extend those of rational
trigonometry in the plane. This gives a unified, computational model of both
spherical and hyperbolic geometries, allows the extension of many results of
Euclidean geometry to the relativistic setting, and provides a new metrical
approach to algebraic geometry.

\end{abstract}

\section*{Introduction}

\textit{Universal Geometry} extends Euclidean and non-Euclidean geometries to
general fields and quadratic forms. This new development is a natural
outgrowth of \textit{Rational Trigonometry }as described in the elementary
text \cite{Wildberger}. It was there developed in the planar context with an
emphasis on the applications to Euclidean geometry. In this paper the subject
is built up in two very general settings---the \textit{affine} one in an $n$
dimensional vector space over a general field with a metrical structure given
by an arbitrary, but fixed, symmetric bilinear form, and the
associated\textit{\ projective} one involving the space of lines through the
origin of a vector space with a symmetric bilinear form.

This allows us to dramatically simplify the usual trigonometric relations for
both Euclidean and non-Euclidean geometries, to extend them to general
bilinear forms, and to reveal the rich geometrical structure of projective
space, with interesting implications for algebraic geometry.

It is pleasant that the main laws of planar rational trigonometry have affine
and projective versions which turn out to hold simultaneously in elliptic
geometry, in hyperbolic geometry, and indeed in any metrical geometry based on
a symmetric bilinear form. The usual dichotomy between spherical and
hyperbolic trigonometry deserves re-evaluation.

The first section introduces and motivates the new approach in a particularly
simple but important special case---that of two dimensional hyperbolic
geometry. The second section establishes the basic trigonometric laws in a
general affine setting, using $n$ dimensional space with an arbitrary
symmetric bilinear form over a general field. The \textit{Triple quad
formula},\textit{\ Pythagoras' theorem}, the \textit{Spread law},
\textit{Thales' theorem}, the \textit{Cross law }and the \textit{Triple spread
formula} include algebraic analogs of the familiar Sine law and Cosine law,
along with the fact that the sum of angles in a triangle is $180^{\circ}$.
Then we derive the projective versions of these laws, which are seen to be
deformations of the affine ones, along with formulas for projective right,
isosceles and equilateral triangles, including general versions of Napier's
rules for solving right triangles in spherical trigonometry. We briefly
mention the important \textit{Spread polynomials}, which are universal analogs
of the Chebyshev polynomials of the first kind, but have an interpretation
over any field. Two examples are shown, one affine and the other projective,
one over the rationals and the other over a finite field. Finally, the Lambert
quadrilateral from hyperbolic geometry is shown to be a general feature in the
projective setting.

Because of the novelty of the approach, some remarks of a subjective nature
may be useful to orient the reader. Metrical geometry is here presented as
fundamentally an \textit{algebraic subject }rather than an \textit{analytic
one}, and the main division in the subject is not between \textit{Euclidean}
and \textit{non-Euclidean}, but rather between \textit{affine} and
\textit{projective}. Elliptic and hyperbolic geometries should be considered
as \textit{projective theories}. Their natural home is the projective space of
a vector space, with metrical structure---\textit{not a metric in the usual
sense}---determined by a bilinear or quadratic form. Over arbitrary fields the
familiar close relation between spheres or hyperboloids and projective space
largely disappears, and the projective space is almost always more basic. The
fundamental formulas and theorems of metrical geometry are those which
\textit{hold over a general field} and are \textit{independent of the choice
of bilinear form}. Many results of Euclidean geometry \textit{extend to the
relativistic setting}, and beyond, once you have understood them in a
universal framework.

This paper lays out the basic tools to begin a dramatic extension of Euclidean
and non-Euclidean geometries, and to once again investigate those aspects of
algebraic geometry concerned with the \textit{metrical properties of curves
and varieties}, in the spirit of Archimedes, Apollonius and the seventeenth
and eighteenth century mathematicians.

\section*{A motivating example from the hyperbolic plane}

Hyperbolic geometry is usually regarded as either a synthetic theory or an
analytic theory. With the synthetic approach you replace Euclid's fifth axiom
with an axiom that allows more than one line through a point parallel to a
given line, and follow Bolyai, Gauss and Lobachevsky as described in
\cite{Greenberg}. With the more modern analytic approach, you introduce a
Riemannian metric in the upper half plane or Poincar\'{e} disk (or sometimes
the hyperboloid of two sheets in three dimensional space), calculate the
geodesics, derive formulas for hyperbolic distances---often employing the
group $PSL_{2}$ of isometries, and then develop hyperbolic trigonometry. This
is described in many places, for example \cite{Beardon} and \cite{Katok}. The
initial interest is in the regular tesselations of the hyperbolic plane,
complex analysis and Riemann surfaces, and the connections with number theory
via quadratic forms and automorphic forms, although nowadays the applications
extend much further.

The new approach to be described here is entirely algebraic and elementary,
and allows us to formulate two dimensional hyperbolic geometry as a
\textit{projective theory} over a general field. There are numerous
computational, pedagogical and conceptual advantages.

Begin with three dimensional space, with typical vector $\left[  x,y,z\right]
$ and bilinear form%
\[
\left[  x_{1},y_{1},z_{1}\right]  \cdot\left[  x_{2},y_{2},z_{2}\right]
=x_{1}x_{2}+y_{1}y_{2}-z_{1}z_{2}.
\]
Define the \textbf{projective point }$a=\left[  x:y:z\right]  $ to be the line
through the origin $O=\left[  0,0,0\right]  $ and the non-zero vector $\left[
x,y,z\right]  .$ The \textbf{projective quadrance} between projective points
$a_{1}=\left[  x_{1}:y_{1}:z_{1}\right]  $ and $a_{2}=\left[  x_{2}%
:y_{2}:z_{2}\right]  $ is the number
\begin{align*}
q\left(  a_{1},a_{2}\right)   & =\frac{\left(  x_{2}y_{1}-x_{1}y_{2}\right)
^{2}-\left(  y_{1}z_{2}-z_{1}y_{2}\right)  ^{2}-\left(  z_{1}x_{2}-x_{1}%
z_{2}\right)  ^{2}}{\left(  x_{1}^{2}+y_{1}^{2}-z_{1}^{2}\right)  \left(
x_{2}^{2}+y_{2}^{2}-z_{2}^{2}\right)  }\\
& =1-\frac{\left(  x_{1}x_{2}+y_{1}y_{2}-z_{1}z_{2}\right)  ^{2}}{\left(
x_{1}^{2}+y_{1}^{2}-z_{1}^{2}\right)  \left(  x_{2}^{2}+y_{2}^{2}-z_{2}%
^{2}\right)  }.
\end{align*}

Define the \textbf{projective line} $L=\left(  l:m:n\right)  $ to be the plane
through the origin (in three dimensional space) with equation%
\[
lx+my-nz=0.
\]
The \textbf{projective spread} between projective lines $L_{1}=\left(
l_{1}:m_{1}:n_{1}\right)  $ and $L_{2}=\left(  l_{2}:m_{2}:n_{2}\right)  $ is
the number
\[
S\left(  L_{1},L_{2}\right)  =\frac{\left(  l_{1}m_{2}-l_{2}m_{1}\right)
^{2}-\left(  m_{1}n_{2}-m_{2}n_{1}\right)  ^{2}-\left(  n_{1}l_{2}-n_{2}%
l_{1}\right)  ^{2}}{\left(  l_{1}^{2}+m_{1}^{2}-n_{1}^{2}\right)
\allowbreak\left(  l_{2}^{2}+m_{2}^{2}-n_{2}^{2}\right)  }.
\]

To give a concrete example, consider the projective triangle $\overline
{a_{1}a_{2}a_{3}}$ with projective points
\[%
\begin{tabular}
[c]{lllll}%
$a_{1}=\left[  1:0:2\right]  $ &  & $a_{2}=\left[  1:-1:3\right]  $ &  &
$a_{3}=\left[  2:1:5\right]  .$%
\end{tabular}
\]
Over the real numbers these lines would correspond to (pairs of) points on the
usual hyperboloid of two sheets inside the null cone $x^{2}+y^{2}-z^{2}=0 $.
The projective lines determined by these projective points are%
\[%
\begin{tabular}
[c]{lllll}%
$L_{1}=a_{2}a_{3}=\left(  8:-1:3\right)  $ &  & $L_{2}=a_{1}a_{3}=\left(
2:1:1\right)  $ &  & $L_{3}=a_{1}a_{2}=\left(  2:-1:1\right)  .$%
\end{tabular}
\]
The projective quadrances of the projective triangle $\overline{a_{1}%
a_{2}a_{3}}$ are then%
\[%
\begin{tabular}
[c]{lllll}%
$q_{1}=-2/5$ &  & $q_{2}=-1/15$ &  & $q_{3}=-4/21$%
\end{tabular}
\]
while the projective spreads are
\[%
\begin{tabular}
[c]{lllll}%
$S_{1}=3/4$ &  & $S_{2}=1/8$ &  & $S_{3}=5/14.$%
\end{tabular}
\]
The analog of the hyperbolic Sine law is
\[
\frac{3/4}{-2/5}=\frac{1/8}{-1/15}=\frac{5/14}{-4/21}%
\]
and there are also analogs of both types of hyperbolic Cosine law.

Over the real numbers, to convert this information into the familiar
Poincar\'{e} model, we intersect the projective point $a=\left[  x_{0}%
:y_{0}:1\right]  $ with the plane $z=1/2,$ yielding the point $A=\left[
x_{0}/2,y_{0}/2,1/2\right]  .$ If $a$ lies inside the null cone $x^{2}%
+y^{2}-z^{2}=0$ then $A$ lies in the open disk $x^{2}+y^{2}<1/4,$ $z=1/2,$
which is the equatorial disk of the sphere $x^{2}+y^{2}+\left(  z-1/2\right)
^{2}=1/4$ with north pole $N=\left[  0,0,1\right]  .$ Project $A$ orthogonally
onto the lower hemisphere of this sphere, giving the point
\[
A^{\prime}=\left[  \frac{x_{0}}{2},\frac{y_{0}}{2},\frac{1}{2}-\frac
{\sqrt{1-x_{0}^{2}-y_{0}^{2}}}{2}\right]
\]
and then stereographically project $A^{\prime}$ from $N$ to the Poincar\'{e}
disk $x^{2}+y^{2}<1,$ $z=0$ (viewed as the open unit disk in the complex
plane) to get the point%
\[
z_{A}=\frac{x_{0}\left(  1-\sqrt{1-x_{0}^{2}-y_{0}^{2}}\right)  }{x_{0}%
^{2}+y_{0}^{2}}+i\frac{y_{0}\left(  1-\sqrt{1-x_{0}^{2}-y_{0}^{2}}\right)
}{x_{0}^{2}+y_{0}^{2}}.
\]
Also project the point $\left[  x_{0},y_{0},1\right]  $ in the plane $z=1$
orthogonally onto the Poincar\'{e} disk to get $w_{A}=\left[  x_{0}%
,y_{0},0\right]  $. Then $z_{A}=\varphi\left(  a\right)  \ $is the
corresponding point in the Poincar\'{e} model to $a$, and $0$ is the
corresponding point to $o=\left[  0:0:1\right]  $. If $\rho$ is the usual
hyperbolic distance in the Poincar\'{e} disk then it is well-known (see for
example \cite[Chapter 7]{Beardon}) that
\[
\rho\left(  0,z_{A}\right)  =\frac{1}{2}\rho\left(  0,w_{A}\right)  =\frac
{1}{2}\log\frac{1+\sqrt{x_{0}^{2}+y_{0}^{2}}}{1-\sqrt{x_{0}^{2}+y_{0}^{2}}}%
\]
and
\[
\sinh^{2}\left(  \rho\left(  0,z_{A}\right)  \right)  =\frac{x_{0}^{2}%
+y_{0}^{2}}{1-x_{0}^{2}-y_{0}^{2}}.
\]
Note that the latter is just the negative of $q\left(  o,a\right)  .$

Now returning to our example, you may use the above formulas to find the
corresponding points in the Poincar\'{e} model to be%
\[%
\begin{tabular}
[c]{lllll}%
$z_{1}=2-\sqrt{3}$ &  & $z_{2}=\frac{3}{2}-\frac{1}{2}\sqrt{7}+i\left(
-\frac{3}{2}+\frac{1}{2}\sqrt{7}\right)  $ &  & $z_{3}=2-\frac{4}{5}\sqrt
{5}+i\left(  1-\frac{2}{5}\sqrt{5}\right)  .$%
\end{tabular}
\]
The standard formula
\[
\rho\left(  z,w\right)  =\log\frac{\left\vert 1-z\overline{w}\right\vert
+\left\vert z-w\right\vert }{\left\vert 1-z\overline{w}\right\vert -\left\vert
z-w\right\vert }%
\]
gives (approximately) the hyperbolic distances
\[%
\begin{tabular}
[c]{lllll}%
$\rho_{1}\approx0.596\,455\,365$ &  & $\rho_{2}\approx0.255\,412\,812$ &  &
$\rho_{3}\approx0.423\,648\,930.$%
\end{tabular}
\]
The corresponding angles in the hyperbolic triangle $\overline{z_{1}z_{2}%
z_{3}}$ may then be calculated using the hyperbolic Cosine Rule%
\[
\cosh\rho_{3}=\cosh\rho_{1}\cosh\rho_{2}-\sinh\rho_{1}\sinh\rho_{2}\cos
\theta_{3}%
\]
to be, in radians, (approximately)%

\[%
\begin{tabular}
[c]{lllll}%
$\theta_{1}\approx2.\,\allowbreak094\,395\,102\approx2\pi/3$ &  & $\theta
_{2}\approx0.361\,367\,126$ &  & $\theta_{3}\approx\allowbreak
0.640\,522\,314.$%
\end{tabular}
\]
To check correctness, you can verify (approximately) the hyperbolic Sine Law%
\[
\frac{\sin2.\,\allowbreak094\,395\,1}{\sinh0.596\,455\,3}\approx\frac
{\sin0.361\,367}{\sinh0.255\,412}\approx\frac{\sin0.640\,522\,}{\sinh
0.423\,648}\approx\allowbreak1.\,\allowbreak36931.
\]
To relate the two approaches, the projective quadrance in the projective
rational model is the \textit{negative} of the square of the hyperbolic sine
of the hyperbolic distance between the corresponding points in the
Poincar\'{e} model, and the projective spread is the square of the sine of the
angle between corresponding geodesics in the Poincar\'{e} model.

The advantages of the projective rational model of the hyperbolic plane
include a cleaner derivation of the theory, simpler and more precise
calculations, with no approximations to transcendental functions required, a
more complete symmetry between rational formulations of the two hyperbolic
Cosine laws, a view of the usual hyperbolic plane as part of a larger picture
involving all projective points, thus accessing the `line at infinity' (the
null cone) as well as the `exterior hyperbolic plane' corresponding to the
hyperboloid of one sheet, and the existence of a beautiful duality between
projective points and projective lines that greatly simplifies hyperbolic
geometry in two dimensions. And as stated previously the theory generalizes to
higher dimensions, to arbitrary fields, and to general symmetric bilinear
forms, and so unifies elliptic and hyperbolic trigonometry.

A more complete account of the two dimensional case, with emphasis on duality,
isometries and applications to tesselations, will be given elsewhere. Now we
turn to develop the general affine theory, and after that the projective
theory which generalizes the above situation.

\section*{Affine rational trigonometry}

Trigonometry studies the \textit{measurement of triangles}. We work in an $n$
dimensional vector space over a field, not of characteristic two. Elements of
the field are referred to as \textbf{numbers}. Elements of the vector space
are called \textbf{points} or \textbf{vectors }(the two terms will be used
interchangeably) and are denoted by $U,V,W$ and so on. The zero vector or
point is denoted $O.$ The unique line $l$ through distinct points $U$ and $V$
is denoted $UV.$ For a non-zero point $U$ the line $OU$ is also denoted
$\left[  U\right]  .$ The unique plane $\Pi$ through non-collinear points
$U,V$ and $W$ is denoted $UVW.$ The plane $OUV$ is also denoted $\left[
U,V\right]  $. The unique affine $3$-flat (translate of a three dimensional
subspace) $\delta$ through non-planar points $U,V,W$ and $Z$ is denoted
$UVWZ,$ and so on.

Fix a symmetric bilinear form and represent it by%
\[
U\cdot V.
\]
In terms of this form, the line $UV$ is \textbf{perpendicular} to the line
$WZ$ precisely when
\[
\left(  V-U\right)  \cdot\left(  Z-W\right)  =0.
\]
A point $U$ is a\textbf{\ null point }or\textbf{\ null vector} precisely when
$U\cdot U=0.$ The \textbf{origin} $O$ is always a null point, but in general
there are others as well. A line $UV$ is a \textbf{null line} precisely when
the vector $U-V$ is a null vector. Null points and lines are themselves well
worth studying, but in this paper we concentrate on non-null points and lines.
Some definitions will be empty when applied to null points or null lines.

For any vector $V$ the number $V\cdot V$ will be denoted $a_{V}$ while more
generally for any vectors $U$ and $V$ the number $U\cdot V$ will be denoted
$b_{UV}.$ Thus $b_{VV}=a_{V}$.

Given vectors $V$ and $U$ the \textbf{projection of }$U$\textbf{\ onto the
line }$\left[  V\right]  $ is
\[
\frac{U\cdot V}{V\cdot V}V.
\]
Then a line perpendicular to $\left[  V\right]  $ in the plane $\left[
U,V\right]  $ is
\[
\left[  U-\frac{U\cdot V}{V\cdot V}V\right]  .
\]

Motivated by this, we define the intersecting planes $\left[  U,W\right]  $
and $\left[  V,W\right]  $ to be \textbf{perpendicular} precisely when
\[
\left(  U-\frac{U\cdot W}{W\cdot W}W\right)  \cdot\left(  V-\frac{V\cdot
W}{W\cdot W}W\right)  =0
\]
which is the same as%
\begin{equation}
\left(  W\cdot W\right)  \left(  U\cdot V\right)  -\left(  U\cdot W\right)
\left(  V\cdot W\right)  =0\label{Perp1}%
\end{equation}
or
\[
a_{W}b_{UV}-b_{UW}b_{VW}=0.
\]

A set $\left\{  U,V,W\right\}  $ of three distinct non-collinear points is a
\textbf{triangle} and is denoted $\overline{UVW}$. The \textbf{lines} of the
triangle $\overline{UVW}$ are $UV,$ $VW$ and $UW.$ Such a triangle is
\textbf{non-null} precisely when each of its lines is non-null.

The \textbf{quadrance} between the points $U$ and $V$ is the number
\[
Q\left(  U,V\right)  =\left(  V-U\right)  \cdot\left(  V-U\right)  .
\]
Then clearly
\begin{align*}
Q\left(  U,V\right)   & =Q\left(  V,U\right)  =V\cdot V-2U\cdot V+U\cdot U\\
& =a_{V}-2b_{UV}+a_{U}.
\end{align*}
The line $UV$ is a null line precisely when $Q\left(  U,V\right)  =0,$ or
equivalently when it is perpendicular to itself.

In ordinary Euclidean geometry distance along a line is additive, assuming you
know what the correct order of points is. With universal geometry, the
important relation between the quadrances formed by three collinear points is
described by a quadratic symmetric function which goes back essentially to
Archimedes' discovery of what is usually called Heron's formula.

\begin{theorem}
[Triple quad formula]If $U,V$ and $W$ are collinear then the quadrances
$Q_{W}=Q\left(  U,V\right)  $, $Q_{U}=Q\left(  V,W\right)  $ and
$Q_{V}=Q\left(  U,W\right)  $ satisfy%
\[
\left(  Q_{U}+Q_{V}+Q_{W}\right)  ^{2}=2\left(  Q_{U}^{2}+Q_{V}^{2}+Q_{W}%
^{2}\right)  .
\]

\end{theorem}

\begin{proof}
First note that the equation can be rewritten as
\[
\left(  Q_{U}+Q_{V}-Q_{W}\right)  ^{2}=4Q_{U}Q_{V}.
\]
Assume $U,V$ and $W$ are collinear and say $U$ and $V$ are distinct so that
$W-U=\lambda\left(  V-U\right)  $ for some number $\lambda.$ In this case%
\begin{align*}
Q_{V}  & =Q\left(  U,W\right)  =\left(  W-U\right)  \cdot\left(  W-U\right)
=\lambda^{2}Q_{W}\\
Q_{U}  & =Q\left(  V,W\right)  =\left(  W-V\right)  \cdot\left(  W-V\right)
=\left(  \lambda-1\right)  ^{2}Q_{W}.
\end{align*}
If $UV$ is a null line then the result is automatic, as both sides are zero,
and otherwise the equation amounts to the identity%
\[
\left(  \left(  \lambda-1\right)  ^{2}+\lambda^{2}-1\right)  ^{2}=4\lambda
^{2}\left(  \lambda-1\right)  ^{2}.
\]

\end{proof}

The next theorem is a restatement and generalization of the most important
theorem in mathematics.

\begin{theorem}
[Pythagoras' theorem]If $U,$ $V$ and $W$ are three distinct points then $UW$
is perpendicular to $VW$ precisely when the quadrances $Q_{W}=Q\left(
U,V\right)  $, $Q_{U}=Q\left(  V,W\right)  $ and $Q_{V}=Q\left(  U,W\right)  $
satisfy%
\[
Q_{W}=Q_{U}+Q_{V}.
\]

\end{theorem}

\begin{proof}
The condition $UW$ perpendicular to $VW$ means that%
\[
\left(  W-U\right)  \cdot\left(  W-V\right)  =0
\]
or%
\[
a_{W}-b_{UW}-b_{VW}+b_{UV}=0.
\]
The condition $Q_{W}=Q_{U}+Q_{V}$ is%
\begin{equation}
a_{V}-2b_{UV}+a_{U}=\left(  a_{V}-2b_{VW}+a_{W}\right)  +\left(  a_{U}%
-2b_{UW}+a_{W}\right)  .\label{Pyth}%
\end{equation}
After a division by two, these conditions are seen to be the same.
\end{proof}

In Euclidean geometry, the separation of lines is traditionally measured by
the transcendental notion of \textit{angle}. The difficulties in defining an
angle precisely, and in extending the concept, are eliminated in rational
trigonometry by using instead the notion of \textit{spread}---in Euclidean
geometry the square of the sine of the angle between two rays lying on those
lines. Fortunately we can formulate the concept completely algebraically.

The \textbf{spread} between the non-null lines $UW$ and $VZ$ is the number%
\[
s\left(  UW,VZ\right)  =1-\frac{\left(  \left(  W-U\right)  \cdot\left(
Z-V\right)  \right)  ^{2}}{Q\left(  U,W\right)  Q\left(  V,Z\right)  }.
\]
This depends only on the two lines, not the choice of points lying on them.
The spread between two non-null lines is $1$ precisely when they are perpendicular.

If $W$ is a non-null point then the \textbf{spread} between the two planes
$OUW$ and $OVW$ may be defined to be the spread between the lines%
\[%
\begin{tabular}
[c]{lllll}%
$\left[  U-\frac{U\cdot W}{W\cdot W}W\right]  $ &  & \textrm{and} &  &
$\left[  V-\frac{V\cdot W}{W\cdot W}W\right]  .$%
\end{tabular}
\]

The next result is an algebraic generalization of the Sine law in planar
trigonometry.$\allowbreak$

\begin{theorem}
[Spread law]Suppose the non-null triangle $\overline{UVW}$ has quadrances
$Q_{W}=Q\left(  U,V\right)  $, $Q_{U}=Q\left(  V,W\right)  $ and
$Q_{V}=Q\left(  U,W\right)  $, and spreads $s_{U}=s\left(  UV,UW\right)  $,
$s_{V}=s\left(  VW,VU\right)  $ and $s_{W}=s\left(  WU,WV\right)  $. Then%
\[
\frac{s_{U}}{Q_{U}}=\frac{s_{V}}{Q_{V}}=\frac{s_{W}}{Q_{W}}.
\]

\end{theorem}

\begin{proof}
Some straightforward simplification shows that$\allowbreak$ the spread
$s_{W}=s\left(  WU,WV\right)  $ is%
\begin{align}
& 1-\frac{\left(  a_{W}-b_{VW}-b_{UW}+b_{UV}\right)  ^{2}}{\left(
a_{W}-2b_{UW}+a_{U}\right)  \left(  a_{W}-2b_{VW}+a_{V}\right)  }%
\label{Num1}\\
& =\frac{\left(
\begin{array}
[c]{c}%
a_{U}a_{V}+a_{U}a_{W}+a_{V}a_{W}+\left(  b_{UV}+b_{UW}+b_{VW}\right)  ^{2}\\
-2\left(  b_{UV}^{2}+b_{UW}^{2}+b_{VW}^{2}\right)  -2a_{U}b_{VW}-2a_{V}%
b_{UW}-2a_{W}b_{UV}%
\end{array}
\right)  }{\left(  a_{W}-2b_{UW}+a_{U}\right)  \left(  a_{W}-2b_{VW}%
+a_{V}\right)  }.\label{Num2}%
\end{align}
The numerator in the expression (\ref{Num2}) is symmetric in the points $U,V$
and $W,$ while the denominator is $Q_{U}Q_{V}.$ It follows that by dividing by
$Q_{W}$ you get an expression which is symmetric in the three points.
\end{proof}

Thales' theorem is the basis of similar triangles in universal geometry.

\begin{theorem}
[Thales' theorem]If the non-null triangle $\overline{UVW}$ has spread
$s_{W}=1$ then
\[
s_{U}=\frac{Q_{U}}{Q_{W}}.
\]

\end{theorem}

\begin{proof}
Immediate from the Spread law.
\end{proof}

The next result is an algebraic generalization of the Cosine law.

\begin{theorem}
[Cross law]Suppose the non-null triangle $\overline{UVW}$ has quadrances
$Q_{W}=Q\left(  U,V\right)  $, $Q_{U}=Q\left(  V,W\right)  $ and
$Q_{V}=Q\left(  U,W\right)  $, and spreads $s_{U}=s\left(  UV,UW\right)  $,
$s_{V}=s\left(  VW,VU\right)  $ and $s_{W}=s\left(  WU,WV\right)  $. Then
\[
\left(  Q_{U}+Q_{V}-Q_{W}\right)  ^{2}=4Q_{U}Q_{V}\left(  1-s_{W}\right)  .
\]

\end{theorem}

\begin{proof}
From (\ref{Pyth}) we have the expression
\begin{align*}
Q_{U}+Q_{V}-Q_{W}  & =\left(  a_{V}-2b_{VW}+a_{W}\right)  +\left(
a_{U}-2b_{UW}+a_{W}\right)  -\left(  a_{V}-2b_{UV}+a_{U}\right) \\
& =2\left(  a_{W}+b_{UV}-b_{UW}-b_{VW}\right)
\end{align*}
while from (\ref{Num1}) we have
\[
1-s_{W}=\frac{\left(  a_{W}-b_{VW}-b_{UW}+b_{UV}\right)  ^{2}}{Q_{U}Q_{V}}.
\]
The result follows.
\end{proof}

Note that the Cross law includes as a special case both the Triple quad
formula and Pythagoras' theorem. The next result is an algebraic analog of the
fact that the sum of the angles in a triangle is $180^{\circ}$.

\begin{theorem}
[Triple spread formula]Suppose the non-null triangle $\overline{UVW}$ has
spreads $s_{U}=s\left(  UV,UW\right)  $, $s_{V}=s\left(  VW,VU\right)  $ and
$s_{W}=s\left(  WU,WV\right)  $. Then%
\[
\left(  s_{U}+s_{V}+s_{W}\right)  ^{2}=2\left(  s_{U}^{2}+s_{V}^{2}+s_{W}%
^{2}\right)  +4s_{U}s_{V}s_{W}.
\]

\end{theorem}

\begin{proof}
Assume the quadrances of $\overline{UVW}$ are $Q_{W}=Q\left(  U,V\right)  $,
$Q_{U}=Q\left(  V,W\right)  $ and $Q_{V}=Q\left(  U,W\right)  $. Write the
Spread law as
\[
\frac{s_{U}}{Q_{U}}=\frac{s_{V}}{Q_{V}}=\frac{s_{W}}{Q_{W}}=\frac{1}{D}%
\]
for some non-zero number $D.$ Now substitute $Q_{U}=Ds_{U},$ $Q_{V}=Ds_{V}$
and $Q_{W}=Ds_{W}$ into the Cross law, and cancel a factor of $D^{2}$,
yielding%
\[
\left(  s_{U}+s_{V}-s_{W}\right)  ^{2}=4s_{U}s_{V}\left(  1-s_{W}\right)  .
\]
Rearrange this to get
\[
\left(  s_{U}+s_{V}+s_{W}\right)  ^{2}=2\left(  s_{U}^{2}+s_{V}^{2}+s_{W}%
^{2}\right)  +4s_{U}s_{V}s_{W}.
\]

\end{proof}

The Triple spread formula can be reinterpreted as a statement about three
non-parallel coplanar lines. If three lines lie in a (two dimensional) plane
then the spread between any two of them is unaffected if one or more of the
lines is translated. In particular we can arrange that the three lines are
concurrent, and so the Triple spread formula still applies.

Another useful observation is that if say $S_{W}=1,$ then the Triple spread
formula becomes%
\[
\left(  s_{U}+s_{V}-1\right)  ^{2}=0
\]
so that
\[
s_{U}+s_{V}=1.
\]

Secondary results of planar rational trigonometry, some developed in
\cite{Wildberger}, are consequences of the main laws of this section, and so
still hold in this general setting.

\section*{An affine example over the rational numbers}

Here is an example of trigonometry in four dimensional space over the rational
numbers (the most important field) with bilinear form
\[
U\cdot V=UMV^{T}%
\]
where
\[
M=%
\begin{pmatrix}
0 & 1 & 0 & 3\\
1 & 1 & 2 & -1\\
0 & 2 & 1 & 0\\
3 & -1 & 0 & -1
\end{pmatrix}
.\allowbreak
\]

Consider the triangle $\overline{UVW}$ where%
\[%
\begin{tabular}
[c]{lllll}%
$U=\left[  1,2,4,3/2\right]  $ &  & $V=\left[  -1,0,1/2,3\right]  $ &  &
$W=\left[  2,2,1,5\right]  .$%
\end{tabular}
\]
Then the quadrances are
\[%
\begin{tabular}
[c]{lllll}%
$Q_{U}=\frac{177}{4}$ &  & $Q_{V}=\frac{71}{4}$ &  & $Q_{W}=\allowbreak38$%
\end{tabular}
\]
and the spreads are%
\[%
\begin{tabular}
[c]{lllll}%
$s_{U}=\frac{10\,263}{10\,792}$ &  & $s_{V}=\allowbreak\frac{3421}{8968}$ &  &
$s_{W}=\frac{3421}{4189}.$%
\end{tabular}
\]
Then you may verify the Spread law%
\[
\frac{10\,263/10\,792}{177/4}=\frac{3421/8968}{71/4}=\frac{3421/4189}%
{38}=\frac{3421}{159182}%
\]
(one of the forms of) the Cross law%
\[
\left(  \frac{177}{4}+\frac{71}{4}-38\right)  ^{2}=4\times\frac{177}{4}%
\times\frac{71}{4}\times\left(  1-\frac{3421}{4189}\right)  .
\]
The Triple spread law becomes%
\begin{align*}
& \left(  \frac{10\,263}{10\,792}+\frac{3421}{8968}+\frac{3421}{4189}\right)
^{2}=\allowbreak\frac{29\,\allowbreak258\,102\,500}{6334\,727\,281}\\
& =2\left(  \left(  \frac{10\,263}{10\,792}\right)  ^{2}+\left(  \frac
{3421}{8968}\right)  ^{2}+\left(  \frac{3421}{4189}\right)  ^{2}\right)
+4\times\frac{10\,263}{10\,792}\times\frac{3421}{8968}\times\frac{3421}{4189}.
\end{align*}

Geometry in such a setting has many familiar features. Here are the
circumcenter $C$ and circumquadrance $K$ of the triangle
\begin{align*}
C  & =\left[  \frac{144}{311},\frac{3789}{3421},\frac{18\,773}{13\,684}%
,\frac{46\,709}{13\,684}\right] \\
K  & =\frac{79\,591}{6842}.
\end{align*}
The orthocenter of the triangle is%
\[
O=\left[  \frac{334}{311},\frac{6106}{3421},\frac{9429}{3421},\frac
{9145}{3421}\right]
\]
the centroid is
\[
G=\left[  \frac{2}{3},\frac{4}{3},\frac{11}{6},\frac{19}{6}\right]
\]
and the nine-point center is%
\[
N=\left[  \frac{239}{311},\frac{9895}{6842},\frac{56\,489}{27\,368}%
,\frac{83\,289}{27\,368}\right]  .
\]
All four points $C,G,N$ and $O$ are collinear, lying on the \textbf{Euler
line} of the triangle, and as expected
\begin{align*}
G  & =\frac{2}{3}C+\frac{1}{3}O\\
N  & =\frac{1}{2}C+\frac{1}{2}O.
\end{align*}
Many other aspects of Euclidean geometry may be explored.

\section*{Projective rational trigonometry}

Fix an $\left(  n+1\right)  $ dimensional vector space over a field with a
symmetric bilinear form $U\cdot V$ as in the previous section. A line through
the origin $O$ will now be called a \textbf{projective point} and denoted by a
small letter such as $u.$ The space of such projective points is called $n$
\textbf{dimensional projective space}. This is a natural domain for algebraic
geometry, and the metrical structure we will introduce gives new directions
for this subject.

If $V$ is a non-zero vector (or point) in the vector space, then $v=\left[
V\right]  $ will denote the projective point $OV.$ A projective point is a
\textbf{null projective point} precisely when some non-zero null point lies on
it. Two projective points $u=\left[  U\right]  $ and $v=\left[  V\right]  $
are\textbf{\ perpendicular} precisely when they are perpendicular as lines.
This is equivalent to the condition $U\cdot V=0.$

A plane through the origin (two dimensional subspace) will be called a
\textbf{projective line }and denoted by a capital letter such as $L.$ A three
dimensional subspace will be called a \textbf{projective plane}. If $V$ and
$W$ are independent vectors then $L=\left[  V,W\right]  $ will denote the
projective line $OVW.$ If $V,W$ and $Z$ are independent vectors then
$\pi=\left[  V,W,Z\right]  $ will denote the projective plane $OVWZ.$

A projective point $u$ \textbf{lies on} a projective line $L$ (or equivalently
$L$ \textbf{passes through} $u$) precisely when the line $u$ lies on the plane
$L$. Similar terminology applies to projective points lying on projective
planes, or projective lines lying on projective planes etc. There is a unique
projective line $L=uv$ which passes through any two distinct projective points
$u$ and $v.$ Three or more projective points which lie on a single projective
line are \textbf{collinear}. Three or more projective lines which all pass
through a single projective point are \textbf{concurrent}.

Two intersecting projective lines $uw$ and $vw$ are \textbf{perpendicular}
precisely when they are perpendicular as intersecting planes. If $u=\left[
U\right]  $, $v=\left[  V\right]  $ and $w=\left[  W\right]  $ then this is
equivalent to the condition (\ref{Perp1})%
\[
\left(  W\cdot W\right)  \left(  U\cdot V\right)  -\left(  U\cdot W\right)
\left(  V\cdot W\right)  =0.
\]

A set $\left\{  u,v,w\right\}  $ of three distinct non-collinear projective
points is a \textbf{projective triangle}, and is denoted $\overline{uvw}.$ The
projective triangle $\overline{uvw}$ is \textbf{null} if one or more of its
projective points is null. The \textbf{projective lines} of the projective
triangle $\overline{uvw}$ are $uv,$ $vw$ and $uw.$

The rich metrical structure of projective space arises with the correct notion
of the separation of two projective points $u$ and $v$. Since each is a line
through the origin in the ambient vector space, we may apply the notion of
spread between these two lines, as developed in the previous section.

The \textbf{projective quadrance} between the non-null projective points
$u=\left[  U\right]  $ and $v=\left[  V\right]  $ is the number
\[
q\left(  u,v\right)  =1-\frac{\left(  U\cdot V\right)  ^{2}}{\left(  U\cdot
U\right)  \left(  V\cdot V\right)  }%
\]
This is the same as the spread $s\left(  OU,OV\right)  $, and has the value
$1$ precisely when the projective points are perpendicular. In terms of the
$a$'s and $b$'s of the previous section
\[
q\left(  u,v\right)  =\frac{a_{U}a_{V}-b_{UV}^{2}}{a_{U}a_{V}}.
\]
Note the use of the small letter $q$ for a projective quadrance to suggest
that it is really a spread, and to distinguish it from a quadrance $Q.$

\begin{theorem}
[Projective triple quad formula]If $u=\left[  U\right]  ,$ $v=\left[
V\right]  $ and $w=\left[  W\right]  $ are collinear projective points then
the projective quadrances $q_{w}=q\left(  u,v\right)  $, $q_{u}=q\left(
v,w\right)  $ and $q_{v}=q\left(  u,w\right)  $ satisfy%
\[
\left(  q_{u}+q_{v}+q_{w}\right)  ^{2}=2\left(  q_{u}^{2}+q_{v}^{2}+q_{w}%
^{2}\right)  +4q_{u}q_{v}q_{w}.
\]

\end{theorem}

\begin{proof}
If $u=\left[  U\right]  ,$ $v=\left[  V\right]  $ and $w=\left[  W\right]  $
are collinear then the three vectors $U,V$ and $W$ are dependent and the lines
$OU,$ $OV$ and $OW$ are coplanar. Since the projective quadrance between
projective points is just the spread between these lines, the theorem is an
immediate consequence of the (affine) Triple spread formula of the previous section.
\end{proof}

Here is the projective version of the most important theorem in mathematics.
Not surprisingly, it has consequences in many directions.

\begin{theorem}
[Projective Pythagoras' theorem]Suppose that $u,v$ and $w$ are three distinct
non-null projective points and that $uw$ is perpendicular to $vw.$ Then the
projective quadrances $q_{w}=q\left(  u,v\right)  $, $q_{u}=q\left(
v,w\right)  $ and $q_{v}=q\left(  u,w\right)  $ satisfy
\[
q_{w}=q_{u}+q_{v}-q_{u}q_{v}.
\]

\end{theorem}

\begin{proof}
If $u=\left[  U\right]  $, $v=\left[  V\right]  $ and $w=\left[  W\right]  $
for some vectors $U,V$ and $W,$ then%
\begin{align*}
& q_{w}-q_{u}-q_{v}+q_{u}q_{v}\\
& =\frac{a_{V}a_{U}-b_{UV}^{2}}{a_{U}a_{V}}-\frac{a_{V}a_{W}-b_{VW}^{2}}%
{a_{V}a_{W}}-\frac{a_{W}a_{U}-b_{UW}^{2}}{a_{U}a_{W}}\\
& +\frac{\left(  a_{V}a_{W}-b_{VW}^{2}\right)  }{a_{V}a_{W}}\frac{\left(
a_{W}a_{U}-b_{UW}^{2}\right)  }{a_{U}a_{W}}\\
& =\frac{\left(  a_{W}b_{UV}+b_{UW}b_{VW}\right)  \allowbreak\left(
b_{UW}b_{VW}-a_{W}b_{UV}\right)  }{a_{W}^{2}a_{V}a_{U}}.
\end{align*}
But we have seen above that the condition that $uw$ is perpendicular to $vw$
is equivalent to%
\[
a_{W}b_{UV}-b_{UW}b_{VW}=0.
\]
Thus this implies that%
\[
q_{w}-q_{u}-q_{v}+q_{u}q_{v}=0.
\]

\end{proof}

Note that the converse does not in general follow.

The \textbf{projective spread} between the intersecting projective lines
$wu=\left[  W,U\right]  $ and $wv=\left[  W,V\right]  $ is defined to be the
spread between these intersecting planes, namely the number
\begin{align*}
& S\left(  wu,wv\right) \\
& =1-\frac{\left(  \left(  U-\frac{U\cdot W}{W\cdot W}W\right)  \cdot\left(
V-\frac{V\cdot W}{W\cdot W}W\right)  \right)  ^{2}}{\left(  \left(
U-\frac{U\cdot W}{W\cdot W}W\right)  \cdot\left(  U-\frac{U\cdot W}{W\cdot
W}W\right)  \right)  \left(  \left(  V-\frac{V\cdot W}{W\cdot W}W\right)
\cdot\left(  V-\frac{V\cdot W}{W\cdot W}W\right)  \right)  }.
\end{align*}
It is undefined if one of the denominators is zero, and is $1$ precisely when
the two projective lines are perpendicular. In terms of the $a$'$s$ and $b$'s%
\begin{equation}
S\left(  wu,wv\right)  =1-\frac{\left(  a_{W}b_{UV}-b_{UW}b_{VW}\right)  ^{2}%
}{\left(  a_{U}a_{W}-b_{UW}^{2}\right)  \left(  a_{V}a_{W}-b_{VW}^{2}\right)
}.\label{SpreadForm2}%
\end{equation}

The projective form of the spread law has the same form as the affine one.

\begin{theorem}
[Projective spread law]Suppose the non-null projective triangle $\overline
{uvw}$ has projective quadrances $q_{u}=q\left(  v,w\right)  $, $q_{v}%
=q\left(  u,w\right)  $ and $q_{w}=q\left(  u,v\right)  $, and projective
spreads $S_{u}=S\left(  uv,uw\right)  $, $S_{v}=S\left(  vw,vu\right)  $ and
$S_{w}=S\left(  wu,wv\right)  $. Then%
\[
\frac{S_{u}}{q_{u}}=\frac{S_{v}}{q_{v}}=\frac{S_{w}}{q_{w}}.
\]

\end{theorem}

\begin{proof}
Assume that $u=\left[  U\right]  ,$ $v=\left[  V\right]  $ and $w=\left[
W\right]  $. After expansion and simplification,%
\begin{equation}
S\left(  wu,wv\right)  =\frac{a_{W}\left(  a_{U}a_{W}a_{V}-a_{V}b_{UW}%
^{2}-a_{U}b_{VW}^{2}-a_{W}b_{UV}^{2}+2b_{UW}b_{UV}b_{VW}\right)  }{\left(
a_{U}a_{W}-b_{UW}^{2}\right)  \left(  a_{V}a_{W}-b_{VW}^{2}\right)
}\label{SpreadForm1}%
\end{equation}
Together with
\[
q_{w}=q\left(  u,v\right)  =\frac{a_{V}a_{U}-b_{UV}^{2}}{a_{U}a_{V}}%
\]
(\ref{SpreadForm1}) shows that the quotient $S_{w}/q_{w}$ is actually
symmetric in the three variables $U,V$ and $W.$
\end{proof}

The next result is a simple but surprising consequence of the Projective
spread law. Even in the simple context of two dimensional elliptic
trigonometry, it reveals that there is an aspect of \textit{similar triangles
in spherical geometry}. This interesting point helps explain why the spread
ratio (opposite quadrance over hypotenuse quadrance) is so important in
rational trigonometry.

\begin{theorem}
[Projective Thales' theorem]If the projective triangle $\overline{uvw}$ has
projective spread $S_{w}=1$ then
\[
S_{u}=\frac{q_{u}}{q_{w}}.
\]

\end{theorem}

\begin{proof}
This follows from the previous theorem.
\end{proof}

Here is the projective version of the Cross law. Unlike the affine result, it
is a \textit{quadratic }equation in the projective spread $S_{w}.$

\begin{theorem}
[Projective cross law]Suppose the non-null projective triangle $\overline
{uvw}$ has projective quadrances $q_{u}=q\left(  v,w\right)  $, $q_{v}%
=q\left(  u,w\right)  $ and $q_{w}=q\left(  u,v\right)  $, and projective
spreads $S_{u}=S\left(  uv,uw\right)  $, $S_{v}=S\left(  vw,vu\right)  $ and
$S_{w}=S\left(  wu,wv\right)  $. Then%
\[
\left(  S_{w}q_{u}q_{v}-q_{u}-q_{v}-q_{w}+2\right)  ^{2}=4\left(
1-q_{u}\right)  \left(  1-q_{v}\right)  \left(  1-q_{w}\right)  .
\]

\end{theorem}

\begin{proof}
If $u=\left[  U\right]  $, $v=\left[  V\right]  $ and $w=\left[  W\right]  $
for some vectors $U,V$ and $W,$ then
\begin{align*}
4\left(  1-q_{u}\right)  \left(  1-q_{v}\right)  \left(  1-q_{w}\right)   &
=4\frac{b_{UV}^{2}}{a_{U}a_{V}}\frac{b_{VW}^{2}}{a_{V}a_{W}}\frac{b_{UW}^{2}%
}{a_{U}a_{W}}\\
& =\frac{4b_{VW}^{2}\allowbreak b_{UW}^{2}\allowbreak b_{UV}^{2}}{a_{W}%
^{2}a_{U}^{2}a_{V}^{2}}.
\end{align*}
Also
\begin{align*}
S_{w}q_{u}q_{v}  & =\frac{a_{W}\left(  a_{U}a_{W}a_{V}-a_{V}b_{UW}^{2}%
-a_{U}b_{VW}^{2}-a_{W}b_{UV}^{2}+2b_{UW}b_{UV}b_{VW}\right)  }{\left(
a_{U}a_{W}-b_{UW}^{2}\right)  \left(  a_{V}a_{W}-b_{VW}^{2}\right)  }\\
& \times\frac{\left(  a_{V}a_{W}-b_{VW}^{2}\right)  }{a_{V}a_{W}}\frac{\left(
a_{U}a_{W}-b_{UW}^{2}\right)  }{a_{U}a_{W}}\\
& =\frac{a_{U}a_{W}a_{V}-a_{V}b_{UW}^{2}-a_{U}b_{VW}^{2}-a_{W}\allowbreak
b_{UV}^{2}+2b_{UW}b_{UV}b_{VW}}{a_{U}a_{V}a_{W}}.
\end{align*}
Thus%
\begin{align*}
& S_{w}q_{u}q_{v}-q_{u}-q_{v}-q_{w}+2\\
& =\frac{a_{U}a_{W}a_{V}-a_{V}b_{UW}^{2}-a_{U}b_{VW}^{2}-a_{W}b_{UV}%
^{2}+2b_{UW}b_{UV}b_{VW}}{a_{U}a_{V}a_{W}}\\
& -\frac{a_{V}a_{W}-b_{VW}^{2}}{a_{V}a_{W}}-\frac{a_{U}a_{W}-b_{UW}^{2}}%
{a_{U}a_{W}}-\frac{a_{V}a_{U}-b_{UV}^{2}}{a_{V}a_{U}}+2\\
& =\frac{2b_{VW}\allowbreak b_{UW}\allowbreak b_{UV}}{a_{W}a_{U}a_{V}}.
\end{align*}
The result follows.
\end{proof}

Note that the \textbf{projective quadrea} $\mathcal{A=}S_{w}q_{u}q_{v}%
=S_{u}q_{v}q_{w}=S_{v}q_{u}q_{w}$ is symmetric in $U,V$ and $W,$ also because
of the Projective spread law.

\begin{theorem}
[Dual projective cross law]Suppose the non-null projective triangle
$\overline{uvw}$ has projective quadrances $q_{u}=q\left(  v,w\right)  $,
$q_{v}=q\left(  u,w\right)  $ and $q_{w}=q\left(  u,v\right)  $, and
projective spreads $S_{u}=S\left(  uv,uw\right)  $, $S_{v}=S\left(
vw,vu\right)  $ and $S_{w}=S\left(  wu,wv\right)  $. Then%
\[
\left(  q_{w}S_{u}S_{v}-S_{u}-S_{v}-S_{w}+2\right)  ^{2}=4\left(
1-S_{u}\right)  \left(  1-S_{v}\right)  \left(  1-S_{w}\right)  .
\]

\end{theorem}

\begin{proof}
If $C_{u}=1-S_{u}$, $C_{v}=1-S_{v}$ and $C_{w}=1-S_{w}$ then the required
identity can be rewritten as%
\begin{equation}
\left(  q_{w}\left(  1-C_{u}\right)  \left(  1-C_{v}\right)  +C_{u}%
+C_{v}+C_{w}-1\right)  ^{2}=4C_{u}C_{v}C_{w}.\label{CrossEq}%
\end{equation}

If $u=\left[  U\right]  $, $v=\left[  V\right]  $ and $w=\left[  W\right]  $
for some vectors $U,V$ and $W,$ then using (\ref{SpreadForm2}) the left hand
side of (\ref{CrossEq}) is the square of an expression which simplifies in a
pleasant fashion to
\[
\frac{2\left(  a_{W}b_{UV}-b_{UW}b_{VW}\right)  \left(  a_{V}b_{UW}%
-b_{UV}b_{VW}\right)  \left(  a_{U}b_{VW}-b_{UV}b_{UW}\right)  }{\left(
a_{V}a_{W}-b_{VW}^{2}\right)  \left(  a_{U}a_{W}-b_{UW}^{2}\right)  \left(
a_{U}a_{V}-b_{UV}^{2}\right)  }.
\]
The square of this is exactly the right hand side of (\ref{CrossEq}).
\end{proof}

\begin{theorem}
[Dual projective Pythagoras' theorem]If $q_{w}=1$ in the notation of the
previous proof, then%
\[
S_{w}=S_{u}+S_{v}-S_{u}S_{v}.
\]

\end{theorem}

\begin{proof}
This follows from the previous result and the polynomial identity%
\begin{align*}
& \left(  S_{u}S_{v}-S_{u}-S_{v}-S_{w}+2\right)  ^{2}-4\left(  1-S_{u}\right)
\left(  1-S_{v}\right)  \left(  1-S_{w}\right) \\
& =\left(  S_{w}-S_{u}-S_{v}+S_{u}S_{v}\right)  ^{2}.
\end{align*}

\end{proof}

The proof of the next result utilized a computer, although it could be checked
by hand.

\begin{theorem}
[Projective triple spread formula]Suppose that $u,v,w$ and $z$ are coplanar
projective points with projective spreads $R_{u}=S\left(  zv,zw\right)  $,
$R_{v}=S\left(  zu,zw\right)  $ and $R_{w}=S\left(  zu,zv\right)  $. Then
\[
\left(  R_{u}+R_{v}+R_{w}\right)  ^{2}=2\left(  R_{u}^{2}+R_{v}^{2}+R_{w}%
^{2}\right)  +4R_{u}R_{v}R_{w}.
\]

\end{theorem}

\begin{proof}
(Using a computer) Suppose that $u=\left[  U\right]  $, $v=\left[  V\right]
$, $w=\left[  W\right]  $ and $z=\left[  Z\right]  $ for vectors $U,V,W$ and
$Z.$ Then%
\begin{align*}
R_{u}  & =\frac{a_{Z}\left(  a_{V}a_{W}a_{Z}-a_{Z}b_{VW}^{2}-a_{V}b_{WZ}%
^{2}-a_{W}b_{VZ}^{2}+2b_{VW}b_{VZ}b_{WZ}\right)  }{\left(  a_{V}a_{Z}%
-b_{VZ}^{2}\right)  \left(  a_{Z}a_{W}-b_{WZ}^{2}\right)  }\\
R_{v}  & =\frac{a_{Z}\left(  a_{U}a_{W}a_{Z}-a_{Z}b_{UW}^{2}-a_{U}b_{WZ}%
^{2}-a_{W}b_{UZ}^{2}+2b_{UW}b_{UZ}b_{WZ}\right)  }{\left(  a_{U}a_{Z}%
-b_{UZ}^{2}\right)  \left(  a_{Z}a_{W}-b_{WZ}^{2}\right)  }\\
R_{w}  & =\frac{a_{Z}\left(  a_{U}a_{V}a_{Z}-a_{Z}b_{UV}^{2}-a_{U}b_{VZ}%
^{2}-a_{V}b_{UZ}^{2}+2b_{UV}b_{UZ}b_{VZ}\right)  }{\left(  a_{U}a_{Z}%
-b_{UZ}^{2}\right)  \left(  a_{Z}a_{V}-b_{VZ}^{2}\right)  }%
\end{align*}
A computer calculation shows that
\[
\left(  R_{u}+R_{v}+R_{w}\right)  ^{2}-2\left(  R_{u}^{2}+R_{v}^{2}+R_{w}%
^{2}\right)  -4R_{u}R_{v}R_{w}%
\]
has a factor which is the determinant
\[%
\begin{vmatrix}
a_{U} & b_{UV} & b_{UW} & b_{UZ}\\
b_{UV} & a_{V} & b_{VW} & b_{VZ}\\
b_{UW} & b_{VW} & a_{W} & b_{WZ}\\
b_{UZ} & b_{VZ} & b_{WZ} & a_{Z}%
\end{vmatrix}
.
\]
But if $u,v,w$ and $z$ are coplanar then $U,V,W$ and $Z$ are linearly
dependent, so this determinant is zero.
\end{proof}

It is worth pointing out that for the case of a non-degenerate bilinear form
and $n=3$ there is a duality between projective points and projective lines,
so the previous three theorems can be deduced from the corresponding earlier
results. However in general this duality is not available.

\section*{Special projective triangles}

The following theorem gives universal analogs of Napier's rules in spherical
trigonometry. The various formulas in this proof are fundamental for
projective trigonometry, and they can all be derived easily enough from the
basic equations of the theorem.

\begin{theorem}
[Napier's rules]Suppose the projective triangle $\overline{uvw}$ has
projective quadrances $q_{u}=q\left(  v,w\right)  $, $q_{v}=q\left(
u,w\right)  $ and $q_{w}=q\left(  u,v\right)  $ and projective spreads
$S_{u}=S\left(  uv,uw\right)  $, $S_{v}=S\left(  vu,vw\right)  $ and
$S_{w}=S\left(  wu,wv\right)  $. If $S_{w}=1$ then any two of the five
quantities $\left[  q_{u},q_{v},q_{w},S_{u},S_{v}\right]  $ determine the
other three, solely through the three \textbf{basic equations}%
\[%
\begin{tabular}
[c]{lllll}%
$q_{w}=q_{u}+q_{v}-q_{u}q_{v}$ &  & $S_{u}=q_{u}/q_{w}$ &  & $S_{v}%
=q_{v}/q_{w}.$%
\end{tabular}
\]

\end{theorem}

\begin{proof}
Two of the projective quadrances allow you to determine the third via the
Projective Pythagoras' theorem $q_{w}=q_{u}+q_{v}-q_{u}q_{v}$, and then the
other two Projective Thales' equations $S_{u}=q_{u}/q_{w}$ and $S_{v}%
=q_{v}/q_{w}$ give the projective spreads.

Given the two projective spreads $S_{u}$ and $S_{v}$, use the Projective
Pythagoras' theorem and the Thales' equations $S_{u}=q_{u}/q_{w}$ and
$S_{v}=q_{v}/q_{w}$ to obtain%
\[
1=S_{u}+S_{v}-S_{u}S_{v}q_{w}.
\]
Thus%
\begin{align*}
q_{u}  & =S_{u}q_{w}=\frac{S_{u}+S_{v}-1}{S_{v}}\\
q_{v}  & =S_{v}q_{w}=\frac{S_{u}+S_{v}-1}{S_{u}}%
\end{align*}
and%
\[
q_{w}=\frac{S_{u}+S_{v}-1}{S_{u}S_{v}}.
\]
If you know a projective spread, say $S_{u},$ and a projective quadrance, then
there are three possibilities. If the projective quadrance is $q_{w},$ then
$q_{u}=S_{u}q_{w}$,%
\[
q_{v}=\frac{q_{w}-q_{u}}{1-q_{u}}=\frac{q_{w}-S_{u}q_{w}}{1-S_{u}q_{w}}%
\]
and
\[
S_{v}=\frac{q_{v}}{q_{w}}=\frac{1-S_{u}}{1-S_{u}q_{w}}.
\]
If the projective quadrance is $q_{u},$ then $q_{w}=q_{u}/S_{u}$,%
\[
q_{v}=\frac{q_{w}-q_{u}}{1-q_{u}}=\frac{q_{u}\left(  1-S_{u}\right)  }%
{S_{u}\left(  1-q_{u}\right)  }%
\]
and%
\[
S_{v}=\frac{q_{v}}{q_{w}}=\frac{1-S_{u}}{1-q_{u}}.
\]
If the projective quadrance is $q_{v},$ then substitute $q_{u}=S_{u}q_{w}$
into the Projective Pythagoras equation to get
\[
q_{w}=S_{u}q_{w}+q_{v}-S_{u}q_{v}q_{w}.
\]
So%
\begin{align*}
q_{w}  & =\frac{q_{v}}{1-S_{u}\left(  1-q_{v}\right)  }\\
q_{u}  & =\frac{S_{u}q_{v}}{1-S_{u}\left(  1-q_{v}\right)  }%
\end{align*}
and
\[
S_{v}=\frac{q_{v}}{q_{w}}=1-S_{u}\left(  1-q_{v}\right)  .
\]

\end{proof}

A projective triangle is \textbf{isosceles} precisely when at least two of its
projective quadrances are equal.

\begin{theorem}
[Pons Asinorum]Suppose a non-null projective triangle $\overline{uvw}$ has
projective quadrances $q_{u}$, $q_{v}$ and $q_{w},$ and projective spreads
$S_{u}$, $S_{v}$ and $S_{w}.$ Then $q_{u}=q_{v}$ precisely when $S_{u}=S_{v}.$
\end{theorem}

\begin{proof}
Follows immediately from the Projective spread law.
\end{proof}

It follows from the Projective Pythagoras' theorem that if $S_{u}=S_{v}=1$
then $q_{w}=0$ or $q_{u}=q_{v}=1.$

\begin{theorem}
[Isosceles projective triangle]Suppose an isosceles projective triangle has
non-zero projective quadrances $q_{u}=q_{v}=q$ and $q_{w},$ and projective
spreads $S_{u}=S_{v}=S$ and $S_{w}.$ Then%
\begin{align*}
q_{w}  & =\frac{4q\left(  1-S\right)  \left(  1-q\right)  }{\left(
1-Sq\right)  ^{2}}\\
S_{w}  & =\frac{4S\left(  1-S\right)  \left(  1-q\right)  }{\left(
1-Sq\right)  ^{2}}.
\end{align*}

\end{theorem}

\begin{proof}
Use the Projective spread law in the form%
\[
S_{w}=\frac{Sq_{w}}{q}%
\]
to replace $S_{w}$ in the Projective cross law%
\[
\left(  q^{2}S_{w}-\left(  2q+q_{w}-2\right)  \right)  ^{2}=\allowbreak
4\left(  1-q\right)  ^{2}\left(  1-q_{w}\right)  .
\]
This yields a quadratic equation in $q_{w}$ with solutions $q_{w}=0,$ which is
impossible by assumption, and
\[
q_{w}=\frac{\allowbreak4q\left(  1-S\right)  \left(  1-q\right)  }{\left(
1-Sq\right)  ^{2}}.
\]
Thus%
\[
S_{w}=\frac{Sq_{w}}{q}=\frac{\allowbreak4S\left(  1-S\right)  \left(
1-q\right)  }{\left(  1-Sq\right)  ^{2}}.
\]

\end{proof}

A projective triangle is \textbf{equilateral }precisely when all its
quadrances are equal. The following formula appeared in the Euclidean
spherical case as Exercise 24.1 in \cite{Wildberger}.

\begin{theorem}
[Equilateral projective triangles]Suppose that a projective triangle is
equilateral with common non-zero projective quadrance $q_{1}=q_{2}=q_{3}=q,$
and with common projective spread $S_{1}=S_{2}=S_{3}=S$. Then%
\[
\left(  1-Sq\right)  ^{2}=4\left(  1-S\right)  \left(  1-q\right)  .
\]

\end{theorem}

\begin{proof}
From the Isosceles projective triangle theorem%
\[
q=\frac{\allowbreak4q\left(  1-S\right)  \left(  1-q\right)  }{\left(
1-Sq\right)  ^{2}}.
\]
Since $q\neq0$ this yields%
\[
\left(  1-Sq\right)  ^{2}=4\left(  1-S\right)  \left(  1-q\right)  .
\]

\end{proof}

The above result is symmetric in $S$ and $q.$ Note that if $S=3/4$ then
$q=8/9.$ This value is important in chemistry---it is the tetrahedral spread
found for example in the methane molecule, and corresponds to an angle which
is approximately $\allowbreak109.\,\allowbreak47^{\circ}.$ As I\ will show
elsewhere, rational trigonometry provides a much more refined analysis of the
geometry of the Platonic solids, but some basic results in this direction can
be found in \cite{Wildberger}.

\section*{Spread polynomials}

We have seen that both affine and projective trigonometry involve the Triple
spread formula%
\[
\left(  a+b+c\right)  ^{2}=2\left(  a^{2}+b^{2}+c^{2}\right)  +4abc.
\]
If $a=b=s$ then $c$ turns out to be either $0$ or $4s\left(  1-s\right)  $. If
$a=4s\left(  1-s\right)  $ and $b=s$ then $c$ turns out to be either $s$ or
$s\left(  3-4s\right)  ^{2}.$ There is then a sequence of polynomials
$S_{n}\left(  s\right)  $ for $n=0,1,2,\cdots$ with the property that
$S_{n-1}\left(  s\right)  ,s$ and $S_{n}\left(  s\right)  $ always satisfy the
Triple spread formula. They play a role in all metrical geometries,
independent of the nature of the symmetric bilinear form, and are defined over
the integers. These are rational analogs of the Chebyshev polynomials of the
first kind, and they have many remarkable properties.

The \textbf{spread polynomial }$S_{n}\left(  s\right)  $ is defined
recursively by $S_{0}\left(  s\right)  =0,S_{1}\left(  s\right)  =s$ and the
rule%
\[
S_{n}\left(  s\right)  =2\left(  1-2s\right)  S_{n-1}\left(  s\right)
-S_{n-2}\left(  s\right)  +2s.
\]

The coefficient of $s^{n}$ in $S_{n}\left(  s\right)  $ is a power of four, so
the degree of the polynomial $S_{n}\left(  s\right)  $ is $n$ in any field of
characteristic not two. It turns out that in the decimal number field
\[
S_{n}\left(  s\right)  =\frac{1-T_{n}\left(  1-2s\right)  }{2}%
\]
where $T_{n}$ is the $n$-th Chebyshev polynomial of the first kind. The first
few spread polynomials are $S_{0}\left(  s\right)  =0$, $S_{1}\left(
s\right)  =s$, $S_{2}\left(  s\right)  =4s-4s^{2}$, $S_{3}\left(  s\right)
=9s-24s^{2}+16s^{3}$, $S_{4}\left(  s\right)  =16s-80s^{2}+128s^{3}-64s^{4}$
and $S_{5}\left(  s\right)  =25s-200s^{2}+560s^{3}-640s^{4}+256s^{5}.$ Note
that $S_{2}\left(  s\right)  $ is the logistic map.

As shown in \cite{Wildberger}, $S_{n}\circ S_{m}=S_{nm}$ for $n,m\geq1,$ and
the spread polynomials have interesting orthogonality properties over finite
fields. S. Goh \cite{Goh} observed that there is a sequence of
`spread-cyclotomic' polynomials $\phi_{k}\left(  s\right)  $ of degree
$\phi\left(  k\right)  $ with integer coefficients such that for any
$n=1,2,3,\cdots$
\[
S_{n}\left(  s\right)  =\prod_{k|n}\phi_{k}\left(  s\right)  .
\]

\section*{A projective example over $\mathbb{F}_{11}$}

Consider a five dimensional vector space over the field $\mathbb{F}_{11}$ with
bilinear form
\[
U\cdot V=UMV^{T}%
\]
where
\[
M=%
\begin{pmatrix}
1 & 10 & 1 & 0 & 0\\
10 & 2 & 5 & 2 & 0\\
1 & 5 & 1 & 4 & 3\\
0 & 2 & 4 & 7 & 2\\
0 & 0 & 3 & 2 & 8
\end{pmatrix}
.
\]
$\allowbreak$ The geometry of the associated four dimensional projective space
is pleasantly accessible. The use of a finite field simplifies calculations
and provides an ideal laboratory for geometrical explorations, even if you are
interested in other fields. Consider the triangle $\overline{uvw}$ where%
\[%
\begin{tabular}
[c]{lllll}%
$u=\left[  1:4:2:6:1\right]  $ &  & $v=\left[  1:2:3:4:1\right]  $ &  &
$w=\left[  0:8:8:3:1\right]  .$%
\end{tabular}
\]

The projective quadrances are then $q_{u}=9$, $q_{v}=8$ and $q_{w}=1,$ and the
projective spreads are $S_{u}=2$, $S_{v}=3$ and $S_{w}=10.$ Since $q_{w}=1,$
$\overline{uvw}$ is a \textbf{dual right triangle},\textbf{\ }so that
$S_{w}=S_{u}+S_{v}-S_{u}S_{v}$.

The Projective spread law is verified to be%
\[
\frac{2}{9}=\frac{3}{8}=\frac{10}{1}%
\]
while the Projective cross law takes the form
\[
\left(  10\times9\times8-9-8-1+2\right)  ^{2}=0=4\left(  1-9\right)  \left(
1-8\right)  \left(  1-1\right)
\]
and the Dual projective cross law takes the form%
\[
\left(  1\times2\times3-2-3-10+2\right)  ^{2}=5=4\left(  1-2\right)  \left(
1-3\right)  \left(  1-10\right)  .
\]
Note that the squares in $F_{11}$ are $0,1,3,4,5,9$, so $S_{v}$ is the only
projective spread of the triangle $\overline{uvw}$ which is a square. As in
the discussion in \cite{Wildberger}, this implies that only the vertex at $v$
has a bisector, and there are two such. The points $b_{1}=\left[
3:0:5:8:7\right]  $ and $b_{2}=\left[  3:2:7:6:10\right]  $ lie on $uw$ and%
\[%
\begin{tabular}
[c]{lllll}%
$S\left(  vu,vb_{1}\right)  =S\left(  vw,vb_{1}\right)  =10$ &  & \textrm{and}
&  & $S\left(  vu,vb_{2}\right)  =S\left(  vw,vb_{2}\right)  =2.$%
\end{tabular}
\]
While bisectors and medians are dependent on number theoretic considerations,
the orthocenter of the projective triangle turns out not to be. You can check
that here it is%
\[
O=\left[  9:1:0:4:1\right]  .
\]

\section*{Lambert quadrilaterals}

Here are two (of many) results from hyperbolic geometry (see \cite[Chapter
7]{Beardon}) that hold more generally.

\begin{theorem}
[Lambert quadrilateral]Suppose the projective points $u,v,w$ and $z$ are
coplanar and form projective spreads
\[
S\left(  uv,uz\right)  =S\left(  vu,vw\right)  =S\left(  wv,wz\right)  =1
\]
and projective quadrances $q\left(  u,v\right)  =q$ and $q\left(  v,w\right)
=p. $ Then%
\[%
\begin{tabular}
[c]{lll}%
$q\left(  w,z\right)  =y=q\left(  1-p\right)  /\left(  1-qp\right)  $ &  &
$q\left(  u,z\right)  =x=p\left(  1-q\right)  /\left(  1-qp\right)  $\\
$q\left(  u,w\right)  =s=q+p-qp$ &  & $q\left(  v,z\right)  =r=\left(
q+p-2qp\right)  /\left(  1-qp\right)  $%
\end{tabular}
\]
and%
\[%
\begin{tabular}
[c]{lllll}%
$S\left(  vu,vz\right)  =x/r$ &  & $S\left(  vw,vz\right)  =y/r$ &  &
$S\left(  wv,wu\right)  =q/s$\\
$S\left(  uw,uv\right)  =p/s$ &  & $S\left(  uw,uz\right)  =q\left(
1-p\right)  /s $ &  & $S\left(  wu,wz\right)  =p\left(  1-q\right)  /s$%
\end{tabular}
\]
and%
\[
S\left(  zu,zw\right)  =S=1-pq.
\]

\end{theorem}

\begin{proof}
The fact that the four points are coplanar implies that the Projective triple
spread theorem applies to any three projective lines of the projective
quadrilateral $\overline{uvwz}$ meeting at a projective point. Furthermore it
implies that where three projective lines meet and one of the spreads is $1,$
the other two spreads must sum to $1.$

The expressions for $S\left(  vu,vz\right)  $, $S\left(  vw,vz\right)  $,
$S\left(  wv,wu\right)  $ and $S\left(  uw,uv\right)  $ follow from the
Projective Thales' theorem. The expression for $s$ follows from the Projective
Pythagoras theorem applied to $\overline{uvw}.$ The same theorem applied to
$\overline{uvz}$ and $\overline{vwz}$ gives the equations%
\begin{align*}
r  & =q+x-qx\\
r  & =p+y-py
\end{align*}
and since $S\left(  vu,vw\right)  =1$
\[
\frac{x}{r}+\frac{y}{r}=1.
\]
These three equations can then be solved to yield the stated values for $r,$
$x$ and $y.$ Also the equations
\[
S\left(  uw,uv\right)  +S\left(  uw,uz\right)  =1=S\left(  wv,wu\right)
+S\left(  wu,wz\right)
\]
can be used to solve for $S\left(  uw,uz\right)  $ and $S\left(  wu,wz\right)
$. Finally use the Projective spread law in $\overline{uwz}$ to get%
\[
S\left(  zu,zw\right)  =S=1-pq.
\]

\end{proof}

\begin{theorem}
[Right hexagon]Suppose a projective hexagon $\overline{a_{1}a_{2}a_{3}%
a_{4}a_{5}a_{6}}$ is planar, meaning that all the projective points lie in
some fixed projective plane, and that all successive projective spreads are
equal to $1,$ that is
\[
S\left(  a_{1}a_{2},a_{1}a_{6}\right)  =S\left(  a_{2}a_{3},a_{2}a_{1}\right)
=\cdots=S\left(  a_{6}a_{1},a_{6}a_{5}\right)  =1.
\]
Then%
\[
\frac{q\left(  a_{1},a_{2}\right)  }{q\left(  a_{4},a_{5}\right)  }%
=\frac{q\left(  a_{2},a_{3}\right)  }{q\left(  a_{5},a_{6}\right)  }%
=\frac{q\left(  a_{3},a_{4}\right)  }{q\left(  a_{1},a_{6}\right)  }.
\]

\end{theorem}

\begin{proof}
Follows by repeated use of the last formula from the previous theorem.
\end{proof}

\end{document}